\documentclass[a4paper, 10pt]{IEEEtran}

\IEEEoverridecommandlockouts                              
\usepackage{amsmath}
\usepackage{cite}

\usepackage[utf8]{inputenc}
\usepackage{latexsym,amsfonts,amssymb,amsmath,amscd,mathrsfs}
\usepackage{graphicx}      
\usepackage{textcomp}
\usepackage[dvipsnames]{xcolor}
\usepackage[shadow,colorinlistoftodos,prependcaption,spanish,textsize=footnotesize]{todonotes}
\usepackage{url}

\newcommand{\vu} {\mathbf{u}}

\newcommand{\setX}{\mathcal{X}}

\newcommand{\setU}{\mathcal{U}}

\newcommand{\R}{\mathcal{R}}

\newtheorem {theo}{Theorem}

\newtheorem {coro}{Corollary}

\newenvironment{prf}{\noindent \textbf{Proof:}}{}

\title{\LARGE \bf Model predictive control for optimal social distancing in a type SIR-switched model}

\author{J. E. Sereno, A. D' Jorge, A. Ferramosca, E.A. Hernandez-Vargas, A. H. Gonz\'alez,
	\thanks{A. H. González, A. D\textquotesingle Jorge and J. E Sereno are with the Institute of Technological Development for the Chemical Industry (INTEC), CONICET-Universidad Nacional del Litoral (UNL), Guemes 3450, Santa Fe (3000), Argentina.}%
	\thanks{A. Ferramosca is with the Department of Management, Information and Production Engineering, University of Bergamo Via Marconi 5, Dalmine (BG) 24044, Italy.}
	\thanks{E.A. Hernandez-Vargas is with the Instituto de Matemáticas, Universidad Nacional Autonoma de Mexico, Boulevard Juriquilla 3001, Santiago de Querétaro, Qro. 76230, Mexico.}%
}

\begin{document}
		\maketitle
	\thispagestyle{empty}
	\pagestyle{empty}

\begin{abstract}                
Social distancing strategies have been adopted by governments to manage the COVID-19 pandemic, since the first outbreak began. However, further epidemic waves keep out the return of economic and social activities to their standard levels of intensity. Social distancing interventions based on control theory are needed to consider a formal dynamic characterization of the implemented SIR-type model to avoid unrealistic objectives and prevent further outbreaks. 
The objective of this work is twofold: to fully understand some dynamical aspects of SIR-type models under control actions (associated with second waves) and, based on it, to propose a switching non-linear model predictive control that optimize the non-pharmaceutical measures strategy. Opposite to other strategies, the objective here is not just to minimize the number of infected individuals at any time, but to minimize the final size of the epidemic while minimizing the time of social restrictions and avoiding the infected prevalence peak to overpass a maximum established by the healthcare system capacity. Simulations illustrate the benefits of the aforementioned proposal.
\end{abstract}


\section{Introduction}

Since its first outbreak in 2019, the novel coronavirus, named as SARS-CoV-2, has paralysed the world social and economic activities. The associated disease, named as COVID-19, widely infected the world population. Considering that initial cases were reported on December 31, 2019, in Wuhan, Chinese province of Hubei, by the date of march 31, 2020, the World Health Organization declared the COVID-19 disease as a pandemic and subsequently, SARS-CoV-2 virus reached the vast majority of countries around the world (\cite{WHO2020}). As a response, governments tackled this first pandemic outbreak by applying social control measures known as Non-Pharmaceutical Interventions (NPIs). These are associated to measures such as social-distancing policy, face-mask requirement, university and school closure, and telework assignment. Certainly, these measures have proved their efficacy to lessens the disease transmission; however, their detrimental effects on social and economic activities have also been showed (\cite{ferguson2020,flaxman2020,ngonghala2020}).

The way the aforementioned measures are decided relies on the so-called SIR-type epidemiological models. SIR-type models are based on the seminal work of \cite{kermack1927}, which firstly established a compartmental relationship between the main variables of an epidemic: Susceptible ($S$), Infected ($I$) and Removed ($C$) individuals, and includes parameters that can be externally modified, as the transmission and recovery/death rates $\beta$ and $\gamma$. Even when several extensions of the original model have been made (to include additional states describing new sub-compartments \cite{giordano2020sidarthe}), the epidemiological objectives always aims to minimize three main indexes (\cite{di2021optimal}): the total fraction infected (or final size of the epidemic, $1-S(\infty)$), the peak of $I$ (or infected peak prevalence) and the average time of the infection.

Both, the minimization of infected peak (\cite{sadeghi2020universal,federico2020taming,morris2021optimal}) and the epidemic final size (\cite{bliman2021best,gonzalez2021dynamic}) were proposed as control objectives in the context of single-interval control strategies; i.e., simple control strategies consisting in a fixed control action that reduces the transmission rate for a finite period of time. One step ahead are the optimal-control-based strategies, that consider a more complex sequence of control actions, together with explicit constraints on both, manipulated and controlled variables. 
In \cite{morato2020optimal} a model predictive controller (MPC) is proposed, based on SIRD (Susceptible-Infected-Recovered-Dead) and SIRASD (SIRD + Asymptomatic-Symptomatic) models. The control objective consists in minimizing both, the number of infected individuals and the time of isolation, while some constraints account for the peak of the infected fraction (maximal values according to the available ICU beds) and for the the dwell-time of control actions (minimal times for both, isolation and no-isolation measurements are imposed to avoid unimplementable policies). 
In \cite{kohler2020robust} - which up to the authors knowledge is the more complete MPC approach - it is proposed a robust non linear MPC based on the SIDARTHE (Susceptible, Infected, Diagnosed, Ailing, Recognized, Threatened, Healed, Extinct) model introduced in \cite{giordano2020sidarthe}. The MPC controller manipulates also the transmission rate, given in this case by several parameters. The objectives are to minimize both, the number of fatalities and the time of isolation, compared to a baseline policy. 
Other similar strategies can be seen in \cite{peni2020nonlinear,carli2020model}. None of the reported works consider the minimization of the final size of the epidemic under finite-time interruption of the measures, in a way that avoid second infection waves, as stated in \cite{gonzalez2021dynamic}.


\subsection{Objective of the work}

The objectives of this manuscript are: (i) to present a dynamical-system perspective to formally analyze SIR-type models, their equilibrium sets and stability, and the second-waves scenarios; and (ii) to propose a Switching Nonlinear Model Predictive Control (swNMPC) strategy that minimizes the epidemic final size, while maintaining the infected peak prevalence under an upper bound imposed by the healthcare system capacity, and minimizing - as long as possible - the social restriction severity.

\section{Review of SIR-type Models}\label{sec:revSIR}

In this work, a non-dimensional version of the typical SIR model (\cite{kermack1927,brauer2012mathematical,sontag2011lecture}) is considered, which can be obtained by rescaling the original time units (\cite{franco2020feedback}): 
\begin{subequations}\label{eq:SIRnondim}
\begin{align}
    \dot{S}(\tau) &= - \R S(\tau) I(\tau) \\
    \dot{I}(\tau) &=   \R S(\tau) I(\tau) - I(\tau) \\
    \dot{C}(\tau) &=  I(\tau),
\end{align}
\end{subequations}
where $S(\tau)$ is the fractions of individuals who are susceptible to contract the infection at time $\tau$, $I(\tau)$ is the fractions of infected individuals (that cause other individuals to become infected), and $C(\tau)$ is the cumulative fractions of removed individuals. $\R$ is the so called \textbf{basic reproduction number}, given by $\R:=\beta/\gamma$, where $\beta$ and $\gamma$ stand for the transmission and the recovery/death rates of the disease, respectively. 
%
%
State variable are constrained to the set
\begin{eqnarray}
\setX := \{(S,I,C) \in \mathbb R^3: S \in [0,1], I \in [0,1], C \in [0,1] \}, \nonumber
\end{eqnarray}
in such a way that $(S(\tau),I(\tau),C(\tau)) \in \setX$ for all $\tau \geq 0$. Furthermore, note that $\dot{S}(\tau)+\dot{I}(\tau)+\dot{C}(\tau) = 0$, and, therefore $S(\tau)+I(\tau)+C(\tau) = 1$, for $\tau \geq 0$. Particularly, $S(0)+I(0)+C(0) = 1$, where $\tau=0$ is assumed to be the epidemic outbreak time, in such a way that $(S(0),I(0),C(0)):=(1-\epsilon,\epsilon,0)$, with $0 < \epsilon \ll 1$, \textit{i.e.}, the fraction of susceptible individuals is smaller than, but close to $1$; the fraction of infected is close to zero and the fraction of removed is null.

The solution of \eqref{eq:SIRnondim} - which was analytically determined in \cite{harko2014exact}, for $\tau \geq \tau_0 >0$, depends on $\R$ and the initial conditions $(S(\tau_0),I(\tau_0),C(\tau_0))\in \setX$. Since $S(\tau) \geq 0$, $I(\tau) \geq 0$, for $\tau \geq \tau_0 >0$, then $S(\tau)$ is a decreasing function of $\tau$ (by (\ref{eq:SIRnondim}.a)) and $C(\tau)$ is an increasing function of $\tau$, for all $\tau \geq \tau_0$.
From (\ref{eq:SIRnondim}.b), it follows that if $S(\tau_0)\R \leq 1$, $\dot{I}(\tau) = (\R S(\tau) -1)I(\tau) \leq 0$ at $\tau_0$. Furthermore, given that $S(\tau)$ is decreasing, $I(\tau)$ is also decreasing for all $\tau \geq \tau_0$.
On the other hand, if $S(\tau_0)\R>1$, $I(\tau)$ initially increases, then reaches a global maximum, and finally decreases to zero. In this latter case, the peak of $I(\tau)$, $\hat I$, is given by $\hat I:=S(\tau_0) + I(\tau_0)-\frac{1}{\R}(1+\ln(S(\tau_0)\R))$, and it is reached at $\hat \tau$, when $\dot I = \R SI - I=0$. This implies the peak of $I$ is reached when $S(\tau)=S(\hat \tau)=S^*$, where 
\begin{eqnarray}\label{eq:Sstar}
S^*:= \min \{1,1/\R\} 
\end{eqnarray}
is a threshold or critical value, known as ``herd immunity". This way, conditions $S(\tau_0)\R>1$ and $S(\tau_0)\R<1$ that determines if $I(\tau)$ increases or decreases at $\tau_0$ can be rewritten as $S(\tau_0)>S^*$ and $S(\tau_0)<S^*$, respectively.

For the sake of simplicity, we define $S_\infty:= \lim_{\tau\rightarrow \infty} S(\tau)$, $I_\infty:= \lim_{\tau \rightarrow \infty} I(\tau)$ and $C_\infty:= \lim_{\tau\rightarrow \infty} C(\tau)$, which are values that depend on initial conditions $S(\tau_0),I(\tau_0)$, $C(\tau_0)$, and $\R$. Taking $\tau \rightarrow \infty$ for the solutions proposed in \cite{harko2014exact}, $I_\infty = 0$, $C_\infty = 1 - S_\infty$ and:
\begin{eqnarray}\label{eq:Ssol2}
	S_\infty(S(\tau_0),I(\tau_0)) = -\frac{W(-\R S(\tau_0)e^{-\R (S(\tau_0)+I(\tau_0))})}{\R}.
\end{eqnarray}
\section{Equilibrium characterization and stability} \label{sec:eq_estabil}

\begin{figure}
	\centering
	\includegraphics[width=0.8\columnwidth]{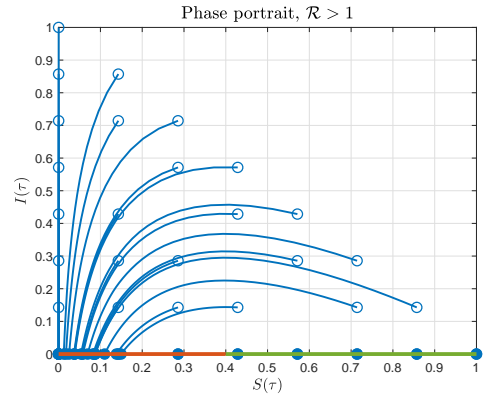}
	\caption{\small{Phase Portrait for System \eqref{eq:SIRnondim} with $\R=2.5$ and starting points summing 1 (starting points in empty circles, ending points in solid circles). Set $\setX_s^{st}$ is in red, while $\setX_s^{un}$ is in green. As it can be seen, all the trajectories converges to $\setX_s^{st}$.}}
	\label{fig:PhaPorRg1}
\end{figure}

\begin{figure}
	\centering
	\includegraphics[width=0.8\columnwidth]{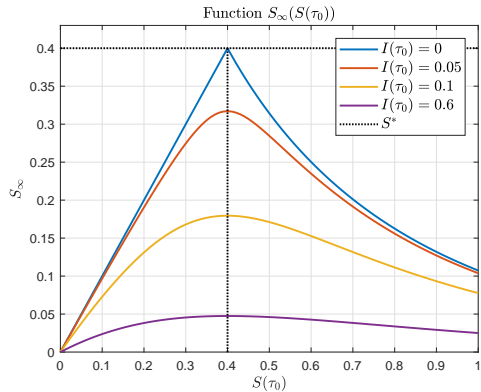}
	\caption{\small{Function $S_\infty(S(\tau_0),I(\tau_0))$, with $\R=2.5$ and fixed values of $I(\tau_0)$.}}
	\label{fig:SinfFunc}
\end{figure}

The equilibria of System \eqref{eq:SIRnondim} is obtained by zeroing each of the differential equations.  For initial conditions $(S(\tau_0),I(\tau_0),C(\tau_0))\in \setX$, this set is given by:
\begin{eqnarray}
\setX_s:=\{(\bar S,\bar I,\bar C) \in \mathbb{R}^3: \bar S \in [0,S(\tau_0)], \bar I=0, \bar C = 1 - \bar S\}. \nonumber
\end{eqnarray}
Next, a key theorem concerning the asymptotic stability of (a subset of) $\setX_s$ is introduced.
\begin{theo}[Asymptotic Stability]\label{theo:stability}
Consider System \eqref{eq:SIRnondim} with arbitrary initial conditions $(S(\tau_0),I(\tau_0),C(\tau_0))\in \setX$, for some $\tau_0 \geq 0$. Then, the minimal (smallest) asymptotically stable equilibrium set in $\setX$ is given by 
\begin{eqnarray}
\setX_s^{st}:=\{(\bar S,\bar I,\bar C) \in \mathbb{R}^3: \bar S \in [0,S^*], \bar I=0, \bar C = 1 - \bar S\}, \nonumber
\end{eqnarray}
while the set
\begin{eqnarray}
\setX_s^{un}:=\{(\bar S,\bar I,\bar C) \in \mathbb{R}^3: \bar S \in (S^*,1], \bar I=0, \bar C = 1 - \bar S\}, \nonumber
\end{eqnarray}
is unstable, being $S^*$ is the herd immunity previously defined.
\end{theo}
\begin{prf}
The proof is given in \cite{gonzalez2021dynamic}.
\end{prf}
Figure~\ref{fig:PhaPorRg1} shows a Phase Portrait for System~\eqref{eq:SIRnondim}, with $\R=2.5$, and initial conditions summing 1.
%
%
%
%
%
\begin{coro}[General system behaviour] \label{cor:sinfty}
Consider System ~\eqref{eq:SIRnondim} with arbitrary initial conditions $(S(\tau_0),I(\tau_0),C(\tau_0))\in \setX$, for some $\tau_0 \geq 0$. Then
\begin{itemize}
    \item $S_\infty \rightarrow 0$ for any value of $S(\tau_0)>0$ and $I(\tau_0)>0$, when $\R \rightarrow \infty$.
    \item $S_\infty$ remains close to $S(\tau_0)$ for any value of $S(\tau_0)>0$ and $I(\tau_0)>0$, when $\R \rightarrow 0$.
    \item If $S(\tau_0) > S^*$ and $I(\tau_0)>0$, $S_\infty$ decreases with $S(\tau_0)$, and $S_\infty < S^*$. 
    \item If $S(\tau_0) < S^*$ and $I(\tau_0)>0$, $S_\infty$ increases with $S(\tau_0)$, and $S_\infty < S^*$.
	\item If $S(\tau_0) = S^*$ and $I(\tau_0) \approx 0$, $S_\infty \approx S^*$, for any value of $\R$ (note that $S^*=1$ for $\R<1$).
\end{itemize}
Figure~\ref{fig:SinfFunc} shows a plot of $S_\infty(S(\tau_0),I(\tau_0))$, corresponding to $\R=2.5$ and fixed values of $I(\tau_0)$.
\end{coro}
%
\section{Social distancing control actions}
Quarantine and isolation, as well as regulations for wearing face masks and avoiding non-essential interactions, are typical non-pharmaceutical measures that the local governments take to reduce infection transmission rates. Social distancing measures affect (decrease) parameter $\beta$ or, directly, parameter $\R$ in System~\eqref{eq:SIRnondim}. In this work, this control action is modeled as follows (\cite{peni2020nonlinear})
\begin{subequations}\label{eq:SIRnondimCont}
	\begin{align}
		\dot{S}(\tau) &= - \R(\tau) S(\tau) I(\tau) \\
		\dot{I}(\tau) &=   \R(\tau) S(\tau) I(\tau) - I(\tau) \\
		\dot{C}(\tau) &=  I(\tau),
	\end{align}
\end{subequations}
where $\R(\tau)$ is a signal ranging from $\R(0)$ (the initial value of $\R$ which is assumed to be $\R(0)>1$) to $\R_{min}$, with $\R(0) > \R_{min} > 0$, which is a minimal value corresponding to the hardest social distancing. In line with other works, it is reasonable to assume that $\R(0) = 3.0$ and $\R_{min}= 0.85$ for simulation proposes (\cite{flaxman2020,li2020,ferguson2020}).


One critical point concerning social distancing measures is that they are always bounded temporary control actions, not permanent ones (as clearly stated in \cite{sadeghi2020universal} and \cite{kohler2020robust}). It is not possible to maintain efficient (or full) social distancing actions for ever (neither for a time long enough to make the infection to disappear) since population fatigue due to psychological or economical problems would systematically relax its effectiveness. 

From an epidemiological point of view, three main general control objectives can be considered (\cite{di2021optimal}): (i) minimize the \textit{epidemic final size (EFS)}: the total fraction infected, $1 - S_\infty$ (\cite{ma2006generality,katriel2012size}), (ii) minimize the \textit{infected peak prevalence (IPP)} or \textit{maximum value of} $I(\tau)$, $\hat I$ (\cite{morris2021optimal,sontag2021explicit}), and (iii) minimize the average \textit{time of infection}, $\bar \tau:= \int_0^\infty \frac{\beta S(t)I(t)}{1 - S(\infty)}dt$.
However, there is not yet a consensus in the optimal control literature about what is better to minimize and to consider as constraints for the system variables. We will assume first a single-interval control action scenario, to reinforce the importance of minimizing the EFS, while keeping the control of the IPP for a second stage. The single-interval control action is given by two times $\tau_i < \tau_f$ defining the initial and final time of social distancing (with $\tau_i<\hat \tau$), and a fixed value of $\R(\tau)=\R_i>0$, to be implemented at $\tau \in [\tau_i,\tau_f]$. For $\tau \in [0,\tau_i) \cap (\tau_f,\infty)$, $\R(\tau)=\R(0)$. In such a context, and according to \cite{gonzalez2021dynamic}, there exist a quasi optimal single interval minimizing the EFS.
\begin{theo}[Quasi optimal single interval control] \label{theo:maxss}
	Consider System~\eqref{eq:SIRnondimCont}, and a given initial time $\tau_i<\hat \tau$. Then, a quasi optimal single interval control (producing $S_\infty \approx S^*$, which is the maximal value $S_\infty$) is given by: $\R_i$ such that $S_\infty(S(\tau_i),I(\tau_i))=S^*$ (we denote this value as $\R_s^{op}$), and $\tau_f$ large enough for the system to approach a QSS ($S(\tau_f) \rightarrow S^*$, $I(\tau_f) \rightarrow 0$).
	Any other single interval control produces $S_\infty<S^*$ and, if $S(\tau_f)>S^*$, a second wave will appear for some $\hat{\hat{\tau}} > \tau_f$.
\end{theo}
\begin{prf}
	The proof is given in \cite{gonzalez2021dynamic}.
\end{prf} \\
\begin{figure}
	\centering
	\includegraphics[width=1\columnwidth]{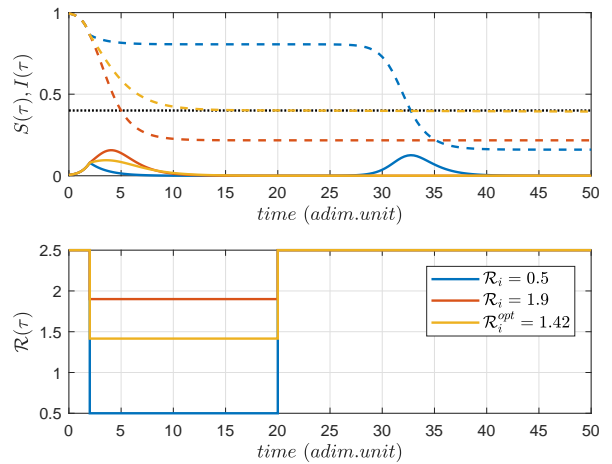}
	\caption{\small{System time evolution corresponding to different single-interval social distancing measures, lasting long enough for the system to reach a QSS before the measures is interrupted. Upper plot: dashed line, $S(\tau)$; solid line, $I(\tau)$; dotted black line, $S^*$. Lower plot: $\R(\tau)$}}
	\label{fig:SingInterSS}
\end{figure}
Figure \ref{fig:SingInterSS} shows the time evolution of $S(\tau)$ and $I(\tau)$ under single-interval social distancing measures lasting long enough for the system to reach a QSS before the measure is interrupted. The initial and final times are given by $\tau_i=2$ and $\tau_f=20$ (non dimensional units), respectively, while $\R(0)=2.5$. For a strong social distancing, $\R_i=0.5 < \R_i^{op}=1.412$ (blue lines), the infected fraction experiences a second wave at $\hat{\hat{\tau}}=32.5$, $12$ times units after the social distancing is interrupted. Furthermore, and more important, $S_\infty = 0.1601$, which is significantly smaller than $S^*=0.4$. For soft social distancing, $\R_i=1.9 > \R_i^{op}=1.42$ (red lines), $S_\infty = 0.2172$, which is also significantly smaller than $S^*$. Finally, for the optimal social distancing, $\R_i^{op}=1.42$, as expected, $S_{\infty}=S^*$ and second waves do not take place. In \cite{gonzalez2021dynamic}, other simulation examples can be seen, concerning short-time single interval controls, i.e., those such that $\tau_f$ is not large enough for the system to reach a QSS. 
As stated in Theorem \ref{theo:maxss}, all these cases produces values of $S_\infty$ significantly smaller than $S^*$. Note that among them, one case is particularly interesting: the one consisting in starting or interrupting any social distancing $\R_i>0$ ($\R=0$ is not a realistic measure) at the very time $\hat \tau$, when $S(\tau)$ reaches $S^*$ from above. This strategy, suggested by some works are far to be optimal from the EFS perspective.

Theorem \ref{theo:maxss} establishes an upper bound for the steady state fraction of susceptible individuals (after the end of social distancing interventions) in the real scenario of temporary control actions of any kind. Furthermore, it states that any social distancing interrupted before a QSS is reached, will produce also a suboptimal value of $S_\infty$. However, the conditions in Theorem \ref{theo:maxss} only determine QSS at the end of the interventions, while the transitory values of the variables remains undetermined. In short, EFS optimal conditions just fixes the area under the curve (AUC) of $I$, but its peaks and general time evolution - corresponding to an arbitrary sequence of social distancing - are not determined. These degrees of freedom for the control actions allows us to propose an optimal control strategy (MPC) that also fulfill maximal conditions for the IPP, $I_{max}$, while minimizing - as long as possible - the side effects of social distancing measures.
\section{Switching NMPC}

\subsection{Switching scheme} \label{sec:switsys}

Social distancing is considered here as a quantified variable that can take only some specific values in a given range. We modeled it as
\begin{eqnarray}\label{eq:Rt}
	\R(\tau) = \R(0) + (\R_{min}-\R(0)) u(\tau), 
\end{eqnarray}
where $u(\tau) \in \setU$ is the control input (manipulated variable), where $\setU=\{0,1/4,1/2,3/4,1\}$ is the set of possible social distancing, being $u=0$ the no social distancing scenario, $u=1$ the hardest social distancing (lockdown), and $u=1/4$, $u=1/2$ and $u=3/4$, intermediary measures. As an example, $u=1/4$ may correspond to social distancing and mask wearing requirements, $u=1/2$ to telework and closure of schools and universities, and $u=3/4$ could be interpreted as a combination of the two. This may have correspondence, for instance, to the classification region system implemented by the Italian Ministry of Health, in which the regions and autonomous provinces are classified into four areas -red, orange, yellow and white- corresponding to three risk scenarios, for which specific restrictive measures are foreseen. 

This way, System~\eqref{eq:SIRnondimCont} together with \eqref{eq:Rt}, can be seen as a switching system, in which the control signal $u$ selects one of 5 possible subsystems, that represents the original system at a given time (\cite{anderson2021discrete,AndersonIFAC2020switching})
\begin{subequations}\label{eq:SIRnondimSig}
	\begin{align}
		\dot{S}(\tau) &= - \R_{\sigma(\tau)} S(\tau) I(\tau) \\
		\dot{I}(\tau) &=   \R_{\sigma(\tau)} S(\tau) I(\tau) - I(\tau) \\
		\dot{C}(\tau) &=  I(\tau),
	\end{align}
\end{subequations}
where $\sigma(\tau)\in \Sigma:=\{1,2,\dots,5\}$ is the switching signal, and $\R_1 = \R(0)$, $\R_2 = \R(0)+\frac{1}{4}(\R_{min}-\R(0))$, $\R_3 = \R(0)+\frac{1}{2}(\R_{min}-\R(0))$, $\R_4 = \R(0)+\frac{3}{4}(\R_{min}-\R(0))$, and $\R_5 = \R_{min}$.

\subsection{Sampled system}\label{sec:sampsys}

To obtain a discrete-time system to be used by the MPC controller, System \eqref{eq:SIRnondimCont} together with \eqref{eq:Rt} is sampled in the general form $x_{k+1} = F(x_{k},u_{k})$, $k \in \mathbb N$, where $F(\cdot)$, is the discrete-time non-linear function, $x_{k}$ is the state vector at the sampled time $k$ (s.t. $x_{k} := x(k T_{s})$, being $T_s>0$ the sampling time), and $u_{k}$ stands for a piece-wise constant input in the form $u(t) = u_{k}, t \in [k T_{s}, (k+1) T_{s}]$. Function $F(\cdot)$ is obtained by the standard explicit Runge–Kutta fourth-order (RK4) method. For the simulation of the non-linear system, the actual trajectory are computed by the backward differentiation formula (BDF) using CVODES solver from the SUNDIALS suite (\cite{hindmarsh2005}).

\subsection{Switching NMPC formulation} \label{sec:swnmpcform}

Here, a model-based controller - which take explicit advantage of the equilibrium/stability characterization (made in Section \ref{sec:eq_estabil}) and the switching scheme (made in Sections \ref{sec:switsys} and \ref{sec:sampsys}) is designed. From an epidemiological perspective, the control objective is to minimize the EFS while maintaining the IPP below an upper bound determined by the health care system and reducing - as much as possible - the time each control action is implemented.
The cost function to be minimized online by the NMPC controller is then given by:
\begin{eqnarray}\label{eq:costorig}
	V_{N} \! \left(x,S^*\!;\! \vu  \right) \!\!=\!\! \sum_{j=0}^{N-1}\!\! Q \|S_j \!\!-\!\! S^* \|^2 \!+\! R\|u_j\|^2 \!+\! P \|S_N \!-\! S^*\|,
\end{eqnarray}
where $N$ is the control horizon, $Q$, $R$ and $P$ are penalizing positive constants, $x=x_k=(S_k,I_k,C_k)$ is the current state at time $k$ and $\vu :=\{u_0,u_1,\cdots,u_{N-1}\}$ is the predicted control sequence.
The optimization problem to be solved at each sampling time $k$ is given by:
\begin{eqnarray*} 
	\begin{array}{rlr}
		\min\limits_{\vu}\;\!\!
		&V_N(x,S^*;\vu) & \\
		s.t.&&   \\
		&x_0=x_k, & \\
		&x_{j+1}= F(x_j,u_j), ~~ \quad \quad j \in \mathbb I_{0:N-1} &  \\
		&x_j \in \setX, ~ \quad \quad \quad \quad \quad \quad \quad  j \in \mathbb I_{0:N-1} &  \\
		&u_j \in \setU, ~ \quad \quad \quad \quad \quad \quad \quad  j \in \mathbb I_{0:N-1} & \\
		&I_j \leq I_{max}, \quad \quad \quad \quad \quad \quad j \in \mathbb I_{0:N-1} &
	\end{array}
\end{eqnarray*}
where $\setU = \{0,0.25,0.5,1\}$ and constraints $I_j \leq I_{max}$ is devoted to impeded that the infected population overpasses a maximum $I_{max}$. Once the optimal solution is computed, then the first optimal input $u_0^0$ is applied to the system, $k \rightarrow k+1$, and the iteration continues with a new solution of the optimization problem (receding horizon control (RHC) strategy).

\begin{figure*}
	\centering
	\includegraphics[width=1 \textwidth]{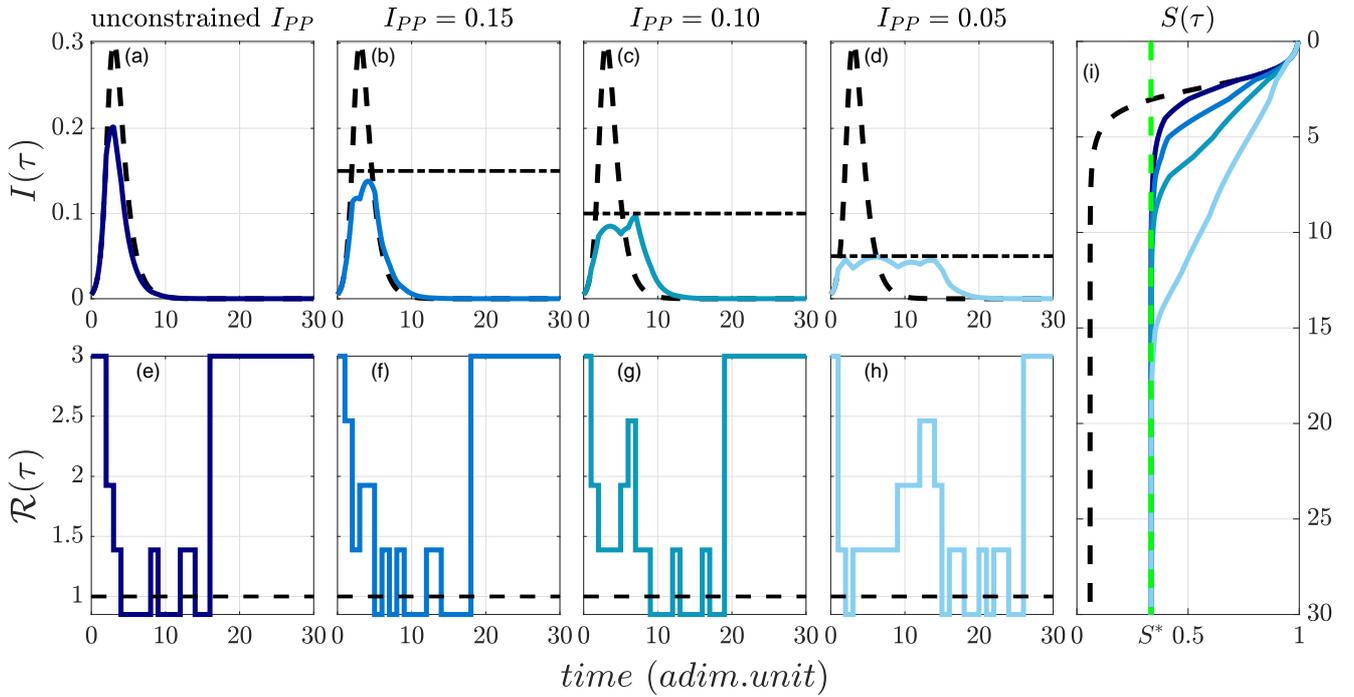}
	\caption{\small{System time evolution corresponding to different restrictions levels of IPP on swNMPC. \textbf{a}-\textbf{d} $I(\tau)$ time evolution under unconstrained IPP (\textbf{a}), $IPP=0.15$  (\textbf{b}), $IPP=0.10$  (\textbf{c}), and $IPP=0.05$  (\textbf{d}). Dash-dotted-black line represent each level of constraint in $IPP$. The switched control signal applied as a social distancing measure $\R_{\sigma(\tau)}$ in each restriction case (\textbf{e}-\textbf{h}), dashed-black horizontal line represent $\R = 1$. From time $\tau=2$ to the end the social distancing intervention is computed by the swNMPC strategy, dashed-black line show the system evolution without interventions. \textbf{i} show how, in each case, the swNMPC reach the optimal herd immunity value $S^*$, dashed black line represent the normalized susceptible population without social distancing measurements. Green-solid vertical line represent the herd immunity value $S^*$.}}
	\label{fig:swMPC}
\end{figure*}

The main advantages of the proposed swNMPC are summarized next. (i) The cost function is null at the optimal behavior (according to the EFS). If other kind of cost functions is used (i.e. to directly minimize $I(\tau)$), this is not necessarily true and the problem will not be well posed. (ii) The upper bound for the IPP, $I_{max}$ can be arbitrary selected.
\section{Simulation Results}

The simulations were run on a laptop computer with i7-4510U (2 cores, 4 threads, 2.0-3.1 GHz) processor and 8GB RAM under MATLAB R2021a using the BONMIN solver, and CasADi version 3.5.5 (\cite{Andersson2019}). We solved the swNMPC optimization Problem~\ref{sec:swnmpcform} under the switched System~\eqref{eq:SIRnondimSig} for a time period $\tau=30$ considering different levels of restrictions for the Infected peak prevalence, $IPP = I_{max}$, even the unconstrained scenario (see Figure~\ref{fig:swMPC} \textbf{a}-\textbf{d}). As it was already said, $\R(0)$ is assumed to be 3, while $\R_{min}=0.85$, The starting time for the social distancing sequence is fixed in 2 times units.

The simulations show that, under the swNMPC social distancing strategy, not further epidemics outbreaks were observed. In fact, this is a key result of the present work. It is important to note that other works present a significant increase of the infected compartments once the control action is dropped and social distancing restrictions interrupted (\cite{peni2020nonlinear,morato2020optimal,carli2020model}). Simulations also show that the harder the infected peak prevalence is, the longer will be the social distancing period. Consequently, unconstrained scenario (Figure~\ref{fig:swMPC}, \textbf{a}-\textbf{e}) has a total of $13$ units of time under social distancing restrictions, and reaches the optimal steady state (herd immunity) at $\tau= 16$. On the other hand, the IPP constraint for $I_{max}= 0.05$, shows a total period of $25$ units of time under social distancing restrictions, and reaches the optimal steady state at $\tau = 25$. Accordingly, these results show that the broader is the capacity of the health care system, the greater is the capacity of governments to manage the pandemic and the shorter is the time spent in lockdown periods. As it can be seen, the strategy tries to minimize - in all cases - the permanence time under strong social distancing ($\R=\R_{min}$).
 
\section{Conclusion}

In this work, a formal dynamical analysis of SIR-type models is presented to consider their equilibrium set and stability into the design of a Switching Model Predictive Control to scheduling social distancing policies. The proposed strategy perform the optimal social distancing policy that minimize the EFS, while maintaining the IPP under an upper bound (according to the healthcare capacity). Discrete levels of social distancing interventions were considered to avoid unrealistic continuous control actions. In the same vein, the IPP maximum were included as a constraint in the NMPC formulation so differences levels of restriction were analyzed showing is possible if the susceptible fraction is under the herd immunity value.

\bibliographystyle{IEEEtran}
\bibliography{bibEpidemio,BSIR_Stabil}

\end{document}